\documentclass[11pt, a4paper]{article}

\usepackage{amsmath,amssymb, amsthm}
\usepackage{mathtools}
\usepackage[utf8]{inputenc}
\usepackage{mathrsfs}
\usepackage{fullpage}
\usepackage{authblk}

\usepackage{dsfont}
\usepackage{bbm}
    
\usepackage{graphicx}
\graphicspath{ {images/} }

\usepackage{wrapfig}

\usepackage{verbatim}

\usepackage{hyperref}

\usepackage{xcolor}
\hypersetup{
  colorlinks   = true, 
  urlcolor     = {blue!90!black}, 
  linkcolor    = {blue!90!black}, 
  citecolor   = {red!90!black} 
}

\newcommand{\R}{\mathbb{R}}
\newcommand{\N}{\mathbb{N}}
\newcommand{\Z}{\mathbb{Z}}
\renewcommand{\d}{\mathrm{d}}

\renewcommand{\P}{\mathbb{P}}

\newcommand{\PP}[1]{\mathbb{P}\left\{#1\right\}}

\newcommand{\Pc}[2]{\mathbb{P}\left[#1 \middle| #2 \right]}

\newcommand{\0}{\mathbf{0}}

\definecolor{ViktorColor}{rgb}{0.17, 0.0, 0.71}

\theoremstyle{plain}
\newtheorem{thm}{Theorem}[section]

\newtheorem{lem}[thm]{Lemma}

\theoremstyle{definition}

\newtheorem{rmk}[thm]{Remark}

\author[1]{Viktor Bezborodov \thanks{Email: \texttt{viktor.bezborodov@pwr.edu.pl}}} 
\author[2]{
Luca Di Persio \thanks{Email: \texttt{luca.dipersio@univr.it}}}
\author[1]{
 Tyll Krueger \thanks{Email: \texttt{tyll.krueger@pwr.wroc.pl}}}

\affil[1]{
	{Wroc\l{}aw University of Science and Technology, Faculty of Electronics  }}
\affil[2]{
	{The University of Verona, Department of Computer Science}}

\title{The continuous-time frog model can spread arbitrarily fast}

\begin{document}

\maketitle

\begin{abstract}
 
The aim of the paper is to demonstrate that 
the continuous-time frog model
can spread arbitrary fast. 
The set of sites visited by an active particle
 can become infinite in a finite time.

\end{abstract}


\section{Introduction}

At time $t = 0$ there are  $\eta (x)$  particles at $x\in \Z ^\d$,
where $\{\eta(x) \}_{x \in \Z ^\d}$ are independent and identically distributed
according to a distribution $\mu$ on $\N \cup \{0\}$.
The particles at the origin are active while all other particles are dormant.
Active particles perform a simple continuous-time random walk
independently of all other particles.
Dormant particles stay put until the first
arrival of an active particle to their site;
upon arrival they become active and start their own simple random walks.
The model was originally defined in discrete time $n = 0,1,2, \dots$
with particles performing a discrete-time simple random walk.
In this paper we consider the continuous-time version.

 In discrete time the frog model cannot spread faster than linearly,
  and the set of locations visited
  by active particles by the time $n \in \N$
  is always contained in  $n\mathcal{D}$, 
  where $\mathcal{D} = \{ (x_1, \dots, x_\d) : |x_1| + \dots + |x_\d| \leq 1 \}$.
In \cite{shapeFrog} and \cite{shapeFrogRandom}
the shape theorem for the discrete-time frog model was established,
and in \cite{shapeFrogRandom}
it was also shown that 
if the tails of $\mu$
are sufficiently heavy, 
 the limiting shape coincides 
 with  $\mathcal{D}$.
 For 
 $\mu = \delta _1$ (delta measure concentrated at $1$)
the shape theorem 
for the continuous-time frog model
was obtained in \cite{stocCombust}.
The frog model has been studied mostly 
in the discrete-time framework.
Recent papers \cite{DHL19} and \cite{benjamini2020}
investigate  respectively 
the coexistence in a two-type
frog model and 
susceptibility properties
on certain finite graphs,
as well as provide an overview of 
other research on this model.  
The transitivity and recurrence properties
of the frog model
attract considerable attention \cite{RecFrog18,   RecFrog17, HJJ16, GNR17}.

In relation to the coexistence 
in two type continuous-time frog model
 the following question was raised
in \cite{DHL19}.
	\begin{flushright}
\begin{minipage}[c]{0,92\textwidth}
	\textbf{Question.} Could the growth be superlinear in time in the continuous time frog
	model if $\eta (x)$ has a very heavy tail?
\end{minipage}
   		\end{flushright}
   	
   In this paper we give a positive answer to this question. Moreover,
   we show that in fact
   for 
   distributions $\mu$ 
   with very heavy tails
   the set of sites visited by active particles
   becomes infinite  in a finite time.
   A precise formulation is given in
   Theorem \ref{short thm}.

\begin{thm}\label{short thm}
	
	There exists a distribution $\mu$
	such that  the time 
	\begin{equation}
	\tau : = \inf \{t: \textnormal{there are infinitely many active particles at } t \}
	\end{equation}
 is a.s. finite.
\end{thm}

We prove Theorem \ref{short thm}
in Section \ref{Sec2}.
In fact,
in Section \ref{Sec2}
we work with  a more general model 
with time between the jumps of random
walks following an arbitrary distribution
rather than the unit exponential.
Remark \ref{rmk waylay} gives an example of an explicit condition on $\mu$ 
ensuring that $\tau < \infty$ a.s.

	The speed of growth
	 of stochastic particle
	systems has been an active field of research 
	for about least half a century
	as the first studies go back at least to the seventies,
	see e.g. \cite{Big76}.
	The superlinear speed for
	a
	branching random walk with
	 polynomial tails
	was demonstrated in \cite{Dur83}.
	The exact speed for the 
	a branching random walk
	satisfying an exponential moment condition
	is given in \cite{Big95, Big97};
	further results 
	and references can be found in \cite{Big_branching_out}.
	More recently  
	the spread rate
	and the maximal displacement
	 of
	modified versions of the model came 
	under investigation.
	A dispersion kernel 
	with tails heavier than 
	exponential but lighter 
	than polynomial
	is treated in 
	\cite{Gan00};
	the spread of the branching random walk
	with certain restrictions 
	is the subject of
	\cite{BM14, trunc_and_crop};
	in
	\cite{FZ12, Mal15}
	the process evolves in a random environment.
	An explosion
		is a phenomenon known
		to take place
		 in  first-passage percolation models
		if a node can have sufficiently many neighbors
		\cite{CD16, ExplFpp}.

		In \cite{frogL}
		conditions ensuring linear or superlinear 
		spread rate
		of the continuous-time frog model
		 are given. 	It turns out
		that whether the spread is  linear or superlinear 
		depends roughly speaking on certain logarithmic moments of $\mu$.			
		The model
		in \cite{Junge20}
		is
		a continuous-time frog model 
		with $m \in \N$ particles per site and  a modified activation 
		mechanism. Specifically, when a site
		belonging to a critical bond percolation cluster 
		is visited for the first time, the 
		sleeping particles (if there are any)
		on the entire cluster are activated. 
	 Thus, many sites can be woken up simultaneously, 
	 and even though there is a fixed number of particles per site, an explosion can occur
	 when the activated clusters are sufficiently large. 
	 In \cite[Theorem 1]{Junge20} the explosion is discussed
	 on $\Z^2 $ and a $d$-ary tree; see also \cite[Theorem 3]{Junge20}.

	In continuous-space settings 
	we  mention a model of growing sets 
	introduced in 
	\cite{Dei03}
	whose speed of growth
	is further studied in \cite{GM08},
	and the spatial birth process \cite{shapenodeath}.
	The linear growth 
	for
	  a discrete-space two-type
	particle model
	is established in
	\cite{KesS05}, see also \cite{KesS08}.
	The model in \cite{KesS05}
	is similar to the frog model,
	however, unlike in the frog model, particles of both types 
	can move. Further discussion takes place in \cite{KRS12}.


\section{The main result, proof, and further discussion} \label{Sec2}

We prove our main result for a generalization
of the frog model in which the particles
perform not  a simple continuous-time random walk,
but a random walk 
with
 the exponential distribution 
of the waiting times between jumps  replaced
by an arbitrary distribution $\pi$ on $(0, \infty)$.
Let $\{ (S _{t} ^{(x ,j)}, t \geq 0 ), x \in \Z ^\d, j \in \N \}$
be
the set of all random walks assigned to particles, $S _{0} ^{(x ,j)} = 0$ for all $x \in \Z ^\d, j \in \N $.
For fixed $t, x$, and $ j$,
$S _{t} ^{(x ,j)} + x$ represents the position 
of  $j$-th particle started at location $x$,
$t$ units of time after the particle was activated.
For each realization of $\eta$, only the walks
$(S _{t} ^{(x ,j)}, t \geq 0 )$
with indices satisfying $j  \leq  \eta (x)$
are used.
For fixed $ x, j$,
the jump times $ \rm{j}_1, \rm{j}_2, \dots$
of $(S _{t} ^{(x ,j)}, t \geq 0 )$
are such that $\rm{j}_{k+1} - \rm{j}_{k}$
are independent random variables  distributed according to $\pi$, $k = 0,1,\dots$ ($\rm{j}_{0} = 0$).
In case of the standard continuous-time frog model, $\pi$ is the unit exponential distribution.

In order
not to exclude  distributions  with an atom at $0$,
we assume that there is at least one active particle at the beginning at the origin $\0 _\d$.
That is,
for realizations of $\eta$ with $\eta (\0 _\d) = 0$ an active particle 
is  added at the origin.

Let us introduce the model
which can serve as a motivation
for treating a more general model rather than only the standard frog model.
Let $\d = 1$. 
Imagine that we again have $ \eta (x)$ particles at $x \in \Z$
at the beginning,
but instead of the random walk the particles now
move
in the continuous space $\R$
 according to independent standard Brownian motions.
Other rules do not change - once some active particles reaches $y\in \Z$
for the first time, all $\eta (y)$ sleeping particles located at $y$ activate
and start their own Brownian motions.
This model can be expressed in 
the discrete-space
 framework
with $\pi$ being the distribution of the time
when the absolute value of a Brownian motion started at $0$
hits $1$, and is thus covered by Theorem \ref{short thm}. Similar models with
particles performing a Brownian motion
were treated in \cite{Brownian_Frog_17, Brownian_Frog_18};
 the description appears already in \cite{stocCombust} in 
 relation to non-isotropy of the lattice models.

Let $\mathcal{A} _t$ be the set of sites visited by an active particle by the time $t$.
If for some  $r >0$, $\pi ( (0, r] ) = 0$, then
for any distribution $\mu$
 a.s. 
$ \mathcal{A} _t \subset [-n, n] ^\d $, where  $n = \lceil \frac{t}{r} \rceil  $,
and hence $(\mathcal{A} _t, t \geq 0)$ grows at most linearly with time.
Lemma \ref{sycophancy} and Theorem \ref{long thm}
show that the reverse is also true.

\begin{lem}  \label{sycophancy}
	Let the dimension $\d = 1$.
	Assume that $\pi ( (0, r] ) > 0$ for all $r >0$.
	  Let $\{A _n\}_{n \in \N}$ and $\{t _n\}_{n \in \N}$ be  increasing sequences of positive numbers,
	$A _n \to \infty$,  $t_n \to t_{\infty} \in (0, \infty]$. 	
	 Then there exists a distribution $\mu$
	such that   $\P \{\sup{ \mathcal{A} _{t_n} }  \geq  A_n \text{ for all } n \in \N \} > 0$,
	and, if $t_\infty < \infty$, the time $\tau$ defined in Theorem \ref{short thm} 
	is a.s. finite.
\end{lem}
The above lemma contains the bulk 
of the proof of Theorem \ref{short thm}. Before proceeding 
to the proof of Lemma \ref{sycophancy} we briefly discuss the main 
idea.
 Take $\{a _n\}_{n \in \N}$
to be a sequence of positive numbers such that
for all $n \in \N$,
\begin{equation}\label{cockup}
\sum _{i=1} ^{n}  a_i \geq A _n.	
\end{equation}
Let  $X_0 = 0$. 
In the proof 
we show that it is possible to construct a sequence of large positive numbers $\{b _n\}_{n \in \N}$
and a distribution $\mu$ such that 
 with positive probability
 there exists a (random) sequence $X_0, X_1, X_2, ...$
 of sites such that the following holds true:

\emph{For every $n \in \N$, there exists a particle started 
at $X_{n-1}$ that moves at least $2 a_{n}$ to the right 
withtin the time $t_n - t_{n-1}$ from its activation, 
and one of the sites in $[X_{n-1} + a_{n},  X_{n-1} + 2 a_{n}]$
contains at least $b_n$ particles. This site is then designated as $X_n$. }

As soon as it is established that such a  sequence $\{X _n\}_{n \in \N}$
exists with positive probability, it can then be deduced
by basically an  ergodicity argument
that with probability one such a sequence \emph{does} exists 
for some possibly different starting site $X_0$.

\textbf{Proof of Lemma \ref{sycophancy}}.
Let $\{a _n\}_{n \in \N}$
be a sequence of positive numbers such that
\eqref{cockup} holds
 for all $n \in \N$ 
and $a_n \geq n^2$.
An example of such a sequence is given by $a_n = A_{n} \vee n^2$.
Let  $t_0 = 0$ and let $\Delta _n = t_{n+1} - t_{n}$ for $n \in \N$.
Define 
$$g(r,m): =  2 ^{-m - 1}  \Bigg( \pi\left( \left(0,  \frac rm \right] \right) \Bigg)^ m, 
\  \  \   r, m > 0,$$
and set $b_n = ( g(\Delta _n, 2 a_n) ) ^{-1} \cdot n $
and let
$\mu$ satisfy
 $\mu([b_{n+1},\infty)) \geq \frac{n}{a_n}$.

Define  a 
random sequence of sites $\{X_n ^{(1)}\}_{n \in \N }$
consecutively as follows:
set $X_0 ^{(1)} = 0$, and
for $n \in \N \cup \{0\}$
set $X_{n+1} ^{(1)} = \infty$
if $X_{n} ^{(1)} = \infty$,
otherwise set 
\begin{multline}
X_{n+1} ^{(1)} = \min\Big\{ k \in \N: a_{n+1} \leq k - X_{n} ^{(1)} \leq 2 a_{n+1} , \eta (k) \geq b_{n+1}, 
\\
\max\{ S _{\Delta _n} ^{(X_{n} ^{(1)} ,j)}: j=1,\ldots, \eta(X_{n} ^{(1)})  \} \geq 2 a_{n+1} \Big\}.
\end{multline}
Here and elsewhere we adopt the convention $\min \varnothing = \infty$.
Let $ \kappa _1 : =  \min\{k \in \N: X_k ^{(1)} = \infty \}  $,
and 
define
\begin{equation} \label{salacioius}
\sigma _1 = 
\begin{cases}
\min\{t \geq t_{\kappa _{_1}} : 
\max\limits _{j = 1,\dots, \eta ( X ^{(1)} _{\kappa _{_1}-1} )}
S _{t  - t_{\kappa _{_1}-1}} ^{(X_{\kappa _{_1}-1} ^{(1)} ,j)} \geq 2 a_{\kappa _{_1}} + 1 \},
&
\text{ on }  \{ \kappa _1 < \infty, \kappa _1 \ne 1 \}, 
\\
\infty, &  \text{ on }  \{ \kappa _1 = \infty \},
\\
\min\{t \geq t_{1} : 
S _{t  } ^{(0 ,1)} \geq 2 a_{1} + 1 \},
&
\text{ on }  \{ \kappa _1 = 1 \}.
\end{cases}   
\end{equation}

Note that a.s.  $\{ \sigma _1 < \infty \} = \{ \kappa _1 < \infty \}$.

We now make the following observation.
If 
the activation 
of some of the sleeping particles
upon coming into contact with an active particle
is delayed or even suppressed entirely,
the resulting process
is going to spread slower
than the frog model.
This also applies to putting to sleep some active particle
and removing (both sleeping and active) particles,
because the spread can only be slowed down as a result.
The slower spread here means
that the set of sites visited by an active particle
by time $t$ for the slowed model
is going to be a subset of the respective set 
for the original model.

Having in mind the  observation above,
we slow down the spread
in multiple ways as described throughout the proof.
The first slowing rule is that at time
$\sigma _0 = 0$
 we 
remove every sleeping particle left of the origin and leave 
a single
 active particle at the origin.
Further,
from time $\sigma _0 $ until  $\sigma _1$
if a site 
with sleeping particles is visited 
by an active particle at time $\theta \in (t_{n-1}, t_{n}]$,
then the sleeping particles at the site
become active and
start moving only after a delay at time $t_{n}$.
Also, 
before time $\sigma _1$
we impose another slowing rule
by restricting the
activation of sleeping particles
to the sites $X_1 ^{(1)}, X_2 ^{(1)}, \ldots$.
Denote by $R_t$ 
the position of the rightmost active particle at time $t$.

On $\{ \kappa _1 < \infty \}$
 at time $\sigma_1$
 we put to sleep every active particle
 keeping only one located at $R_{\sigma_1}$,
 and restart the process
 in the same fashion.
 (We note here that given $\{\sigma _1 < \infty \}$, the random variables 
 $\eta(R_{\sigma _1} + 1), \eta(R_{\sigma _1} + 2), \dots $
 are  independent and distributed according to $\mu$.
 Thus, the usage of the word `restart' is justified 
 as the restarted process is going to have the same distribution.)

 Define the sequence 
 $\{X_n ^{(2)}\}_{n \in \N }$
 by
 setting  $X_0 ^{(2)} = R_{\sigma  _1}$
  on the event $\{ \sigma _1 < \infty \} $
  and $X_0 ^{(2)} = \infty$ 
  on the complement  $\{ \sigma _1 < \infty \} ^{c} = \{ \sigma _1 = \infty \}  $, and
 for $n \in \N \cup \{ 0\}$
 by setting $X_{n+1} ^{(2)} = \infty$
 if $X_{n} ^{(2)} = \infty$,
 and
 otherwise 
 \begin{multline}
 X_{n+1} ^{(2)} = \min\big\{ k \in \N: a_{n+1} \leq k - X_{n} ^{(2)} \leq 2 a_{n+1} , \eta (k) \geq b_{n+1}, 
 \\
 \max\{ S _{\Delta _n} ^{(X_{n} ^{(2)} ,j)}: j=1,\ldots, \eta(X_{n} ^{(2)})  \} \geq 2 a_{n+1}    \big\}.
 \end{multline}
 We then define
 $ \kappa _2 : =  \min\{k \in \N \cup \{0 \}: X_k ^{(1)} = \infty \}  $
 and set
 \begin{equation} \label{salacioius2}
 \sigma _2 = 
 \begin{cases}
 \min\{t \geq t_{\kappa _{_2}} +  \sigma _{1} : 
 \max\limits _{j = 1,\dots, \eta ( X ^{(2)} _{\kappa _{_2}-1} )}
 S _{t - \sigma _{1} - t_{\kappa _{_2}-1}} ^{(X_{\kappa _{_2}-1} ^{(2)} ,j)} \geq  2 a_{\kappa _{_2}} + 1 \} 
 &
 \text{ on }  \{1 < \kappa _2 < \infty \}, 
 \\
 \infty, &  \text{ on }  \{ \kappa _2 = \infty \},
 \\
 1 \ (\text{this value is arbitrary and does not affect anything}), &  \text{ on }  \{ \kappa _2 = 0 \},
 \\
 \min\{t \geq t_{1} +  \sigma _{1} : 
 S _{t - \sigma _{1} } ^{(X_{0} ^{(2)} ,1)} \geq 2 a_{1} + 1 \},
 &
 \text{ on }  \{ \kappa _2 = 1 \}.
 \end{cases}   
 \end{equation}
 Next
define the sequences
 $\{X_n ^{(3)}\}_{n \in \N }$,
 $\{X_n ^{(4)}\}_{n \in \N }$, $\dots$,
 and the times $\kappa _3$, $\sigma _3$, $\dots$,
  consecutively in the same fashion.
 
  On  $\{ \sigma _1 < \infty \} $,
  the same restrictions are introduced on
  the time interval
   $(\sigma _1, \sigma _2]$
  as on $(\sigma _0, \sigma _1]$.
  Specifically,
 at $\sigma _1$ every sleeping particle left to $X_0 ^{(2)} = R_{\sigma _1}$
 is removed.
  From time $\sigma _1$ until  $\sigma _2$,
 the activation of sleeping particles
 at a site first visited by an active particles 
 during the time interval $(\sigma _1 + t_{n-1}, \sigma _1 + t_{n}]$
 takes place with a delay  at  $\sigma _1 + t_n$.
 The activation of sleeping particles is only allowed on sites
 $X_1 ^{(2)}, X_2 ^{(2)}, \ldots$.
On  $\{ \sigma _1 < \infty \} \cap \{ \sigma _2 < \infty \} $,
 same restrictions are made during $(\sigma _2, \sigma _3]$, and so on.


For $n, m \in \N$
denote  $Q_n ^{(m)} = \{X_n ^{(m)} < \infty \} $, 
and let $Q_ \infty ^{(m)}  = \bigcap\limits _{n \in \N} Q_n ^{(m)} = \lim\limits _{n \to \infty} Q_n ^{(m)} $
be the event $\{X_k ^{(m)} < \infty, k \in \N \} = \{\kappa _m = \infty \}$ that 
all elements of the  sequence $\{X_n ^{(m)}\}_{n \in \N }$ are finite.
By construction $\eta (X_n ^{(1)}) \geq b_n$ a.s. on $Q_n ^{(1)}$, 
hence 
by Lemma \ref{lapidation}
\begin{multline}\label{tiff}
\Pc{ \max\{ S _{\Delta _n} ^{(X_{n} ^{(1)} ,j)}: j=1,\ldots,
	 \eta(X_{n} ^{(1)})  \} \geq 2 a_{n+1}  }{Q_n ^{(1)}} \geq 
  1 - \left[ 1 - \PP{S_{\Delta _n} \geq 2a_n}  \right] ^{b_n}
  \\
  \geq 1 - \left[1 -
  g (\Delta_n, 2 a_n)
    \right] ^{b_n} 
  \geq 1 - \left[1 -  g (\Delta_n, 2 a_n) \right] 
  ^{ (g (\Delta_n, 2 a_n)) ^{-1} \cdot n }
  \geq 1 - e^{-n}.
\end{multline}
In \eqref{tiff} we used the inequality $\left(1 - \frac 1y \right) ^y < e ^{-1}$ for $y > 1$.
At the same time we have
\begin{multline}\label{delapidated}
\Pc{ \eta (X_n ^{(1)} + a_{n+1}) \vee  \eta (X_n ^{(1)} + a_{n+2}) \vee \ldots   \vee  \eta (X_n ^{(1)} + 2 a_{n + 1}) \geq b_{n+1} }
{Q_n ^{(1)}} 
\\
\geq  1 - \left[ 1 - \mu ([b_{n+1}, \infty ) )  \right] ^{a_n} 
\geq  1 - \left[ 1 - \frac{n}{a_n}  \right] ^{a_n} 
\geq 1 - e^{-n}.
\end{multline}
Since 
\begin{align}
 Q_{n+1} ^{(1)} =\ & Q_n ^{(1)} \cap \{\max\{ S _{\Delta _n} ^{(X_{n} ^{(1)} ,j)}: j=1,\ldots,
 \eta(X_{n} ^{(1)})  \} \geq 2 a_{n+1}  \}
 \\
 & \cap \{ \eta (X_n ^{(1)} + a_{n+1} ^{(1)}) \vee  \eta (X_n ^{(1)} + a_{n+2} ^{(1)}) \vee \ldots   \vee  \eta (X_n ^{(1)} + 2 a_{n + 1}) \geq b_{n+1}\}, 
 \notag
\end{align}
by \eqref{tiff} and \eqref{delapidated}
\begin{equation}
\Pc{ Q_{n+1} ^{(1)} }{Q_n ^{(1)} } \geq 1 - 2e^{-n}.
\end{equation}
Hence
\begin{equation} \label{trepidation}
\PP{Q _{\infty} ^{(1)} } = \lim\limits _{n \to \infty} \PP{Q_n ^{(1)}} = 
\PP{Q_1 ^{(1)}} \prod\limits _{n = 1} ^\infty \Pc{Q_{n+1} ^{(1)}}{Q_n ^{(1)}} \geq \PP{Q_1 ^{(1)}} \prod\limits _{n = 1} ^\infty ( 1 - 2e^{-n}) >0.
\end{equation}
A.s.
on $Q _{\infty} ^{(1)} $,
 $\sup{\mathcal{A} _{t_n}} \geq X_n \geq \sum\limits_{i= 1} ^n a_i \geq A_n$,
 so the first statement of the lemma is proven.

Let $Q  ^{\infty} = \bigcup \limits _{m= 1} ^\infty (Q _\infty ^{(m)})
 =\{ \kappa _m = \infty \text{ for some } m\in \N \} $
 be the event that for some $m\in \N$,
 all elements of the sequence $\{X_n ^{(m)}\}_{n \in \N }$ are finite.
Now we can use a standard restart argument to show that
$\PP{ (Q  ^{\infty}) ^{c} } = 0$, that is $\PP{ Q  ^{\infty} } = 1$.
Because of the independence of the random walks,
the distribution of 
$\{X_n ^{(m+1)} - R_{\sigma _{m}}  \}_{n \in \N }$
given 
$ \bigcap\limits _{i = 1} ^m(Q _\infty ^{(i)})^{c}$
coincides with the (unconditional) distribution of 
${\{X_n ^{(1)} - R_{\sigma _{0}} \}_{n \in \N } = \{X_n ^{(1)}  \}_{n \in \N }}$.
Hence by \eqref{trepidation}
\begin{equation}
\begin{aligned}
 \PP{ (Q  ^{\infty}) ^{c} } 
 = \PP{ \bigcap \limits _{m= 1} ^\infty (Q _\infty ^{(m)})^{c}  }
   = & \ \PP{(Q _\infty ^{(1)})^{c}}  \prod\limits _{m = 1} ^\infty \Pc{(Q _\infty ^{(m + 1)})^{c}  }
   { \bigcap\limits _{i = 1} ^m(Q _\infty ^{(i)})^{c} }
   \\
   =  & \ \PP{(Q _\infty ^{(1)})^{c}} \prod\limits _{m = 1} ^\infty \left[ 1 - \PP{Q _{\infty} ^{(1)} } \right]
    = 0.
   \end{aligned}
\end{equation}

Thus $\PP{ (Q  ^{\infty})  } = 1$, consequently a.s. there exists $m \in \N$
such that
 the elements of the sequence $\{X_n ^{(m)}\}_{n \in \N }$
are
all
 finite  and $ \kappa _m = \infty$. 
 Note that  this implies that a.s. $\kappa _1, \dots , \kappa _{m-1} < \infty$ if $m>1$.
 In particular, a.s. on $\{ m >1 \}$ we have $\sigma _{m-1} < \infty$. 
By 
construction
the sites $X_1 ^{(m)}$, $X_2 ^{(m)}$, $\dots$,
 are occupied 
 at the time $\sigma _{m-1} + t _1$, $\sigma _{m-1} + t _2$, $\dots$,
 respectively,
 and 
 $X_{n+1} ^{(m)} - X_n ^{(m)} \geq a_{n+1}$, $n \in \N \cup\{0\}$.
 Thus an infinite number of sites 
 have been visited by an active particle by the  time $\sigma _{m-1} + t _\infty$,
 which is a.s. finite if $t_\infty < \infty$.
 \qed

\begin{thm} \label{long thm}
	Assume that $\pi ( (0, r] ) > 0$ for all $r >0$. Then 
 		there exists a distribution $\mu$
	such that  the time 
 $\tau$ defined in Theorem \ref{short thm}
	is a.s. finite.
\end{thm}
\textbf{Proof}.
The one-dimensional projections 
of the particles of the $\d$-dimensional model
perform a random walk
whose times between jumps
are
distributed according to
$\pi ^{(1)}  = \sum\limits _{n=1} ^ \infty \frac1 \d (\frac{\d-1}{\d})^{n-1} \pi ^{\ast n} $.
Hence
the projection of a continuous-time $\d$-dimensional frog model
on an axis
dominates
a continuous-time one-dimensional frog model
having $\pi ^{(1)}  $
as the distribution between jumps of random walks
and the same initial sleeping particles distribution $\mu$. 

Specifically,
 recall that $\mathcal{A}_t$ 
is the set of sites visited by an active particle
by time $t$ for the $\d$-dimensional
 frog models with 
 time intervals between jumps distributed according to $\pi$, 
 and let 
 $\mathcal{A}_t ^{(1)}$
be the sets of sites visited by an active particle
by time $t$ for the
one-dimensional frog models
with 
 intervals between jumps distributed according to $\pi ^{(1)} $.
 Then
 $(\mathcal{A}_t, t\geq 0 )$ 
 and 
 $( \mathcal{A}_t ^{(1)}, t\geq 0 )$
 can be coupled in such a way that
  a.s. $\Pi _1 \mathcal{A}_t \supset \mathcal{A}_t ^{(1)}$
  for all
$t \geq 0$,
where $\Pi _1$ is the projection on the first coordinate.
Since $\pi ^{(1)}  $
satisfies conditions of Lemma \ref{sycophancy}
if  $\pi$ satisfies conditions of Theorem \ref{long thm},
by Lemma \ref{sycophancy}
the set
$\mathcal{A}_s ^{(1)}$
is infinite for some $s \in (0,\infty)$.
Hence so is $\mathcal{A}_s $.
\qed

Theorem \ref{short thm} for the standard frog model is a  particular case of Theorem \ref{long thm}.

\begin{rmk}
	If the dimension $\d = 1$,
	then $\tau < \infty$ a.s. 
	implies by symmetry  that 
	the time when every site has been visited
	by an active particle, i.e. the time
	there are no sleeping particles left, is also a.s. finite.
	In the terminology of \cite{benjamini2020}
	the model is susceptible, despite the underlying graph being infinite.
	Extending this to higher dimensions 
	and other graphs
	 would require 
	additional arguments.
\end{rmk}

\begin{rmk}
	It follows from
	the proof of
	 Lemma \ref{sycophancy}
	that for every $\varepsilon > 0$,
	$\mu$ can be chosen in such a way that 
	\begin{equation}
	 \PP{ \tau > \varepsilon} \leq  \varepsilon.
	\end{equation}
\end{rmk}

\begin{rmk}\label{rmk waylay}
	Taking $\pi$ to be the unit exponential distribution, $a_n = n^2$,
	$\Delta_n = \frac{1}{n^2}$,
	and $b_n = 2^{4n^2 + 1}n ^{8n^2 + 1}$, we see that the conditions
	in the proof 
	of Lemma 
   \ref{sycophancy}
   are satisfied and $t_\infty <  \infty$. Thus, an example of an explicit condition  on $\mu$ 
   implying $\tau < \infty$ 
   is given by $\mu( [2^{4n^2 + 1}n ^{8n^2 + 1}, \infty)) \geq \frac{1}{n-1} $,
   $n\geq 2$.
\end{rmk}


The next lemma provides a lower estimate
of the tails of a random walk
performed by an active particle.
It is used in the proof of Lemma \ref{sycophancy}.

\begin{lem} \label{lapidation}
	Let $(S_t, t\geq 0 )$ be a  continuous-time random walk
	on $\Z$, $S_0 = 0$, with times between jumps distributed according to $\pi$,
	 and let $r >0$. Then 
	\begin{equation}
	\PP{S_r \geq m} \geq 2 ^{-m - 1}  \Bigg( \pi\left( \left(0,  \frac rm \right] \right) \Bigg)^ m.
	\end{equation}
\end{lem}
 \textbf{Proof}.
 Let $\rm{j}_k$ be the time of the $k$-th jump of $(S_t, t\geq 0 )$.
  Since the direction and the timing of each jump are independent,
 \begin{align*}
  \PP{S_r \geq m} \geq & \ 
   \P
 \Bigg\{ \rm{j}_1 \leq \frac rm, \rm{j}_2 -\rm{j}_1 \leq \frac rm,  \dots, 
  \rm{j}_m -\rm{j}_{m-1} \leq \frac rm
	\Bigg\}
\\
&
\times 
\PP{ \text{first } m \text{ jumps are all to the right} 
} \times \PP{ S_r - S _{\rm{j}_m} \geq 0  } 
\\ 
&
\geq \Bigg( \pi\left( \left(0,  \frac rm \right] \right) \Bigg)^ m   2 ^{-m} \frac 12
=
 2 ^{-m - 1}  \Bigg( \pi\left( \left(0,  \frac rm \right] \right) \Bigg)^ m.
  \qedhere
 \end{align*}
\qed

\section*{Acknowledgements}

Viktor Bezborodov is grateful for the support 
from the University of Verona.

\bibliographystyle{alphaSinus}
\bibliography{Sinus}

\newcommand{\etalchar}[1]{$^{#1}$}
\begin{thebibliography}{BDPK{\etalchar{+}}17}

\bibitem[AMP02]{shapeFrog}
O.~S.~M. {Alves}, F.~P. {Machado}, and S.~Y. {Popov}.
\newblock {The shape theorem for the frog model.}
\newblock {\em {Ann. Appl. Probab.}}, 12(2):533--546, 2002.

\bibitem[AMPR01]{shapeFrogRandom}
O.~S.~M. {Alves}, F.~P. {Machado}, S.~Y. {Popov}, and K. {Ravishankar}.
\newblock {The shape theorem for the frog model with random initial
  configuration.}
\newblock {\em {Markov Process. Relat. Fields}}, 7(4):525--539, 2001.

\bibitem[BDD{\etalchar{+}}18]{Brownian_Frog_18}
E. {Beckman}, E. {Dinan}, R. {Durrett}, R. {Huo}, and M. {Junge}.
\newblock {Asymptotic behavior of the Brownian frog model.}
\newblock {\em {Electron. J. Probab.}}, 23:19, 2018.
\newblock Id/No 104.

\bibitem[BDPK{\etalchar{+}}17]{shapenodeath}
V. Bezborodov, L. Di~Persio, T. Krueger, M. Lebid, and T. O\.za\'nski.
\newblock Asymptotic shape and the speed of propagation of continuous-time
  continuous-space birth processes.
\newblock {\em Advances in Applied Probability}, 50(1):74–101, 2017.

\bibitem[BDPKT20]{trunc_and_crop}
V. Bezborodov, L. Di~Persio, T. Krueger, and P. Tkachov.
\newblock Spatial growth processes with long range dispersion: Microscopics,
  mesoscopics and discrepancy in spread rate.
\newblock {\em Ann. Appl. Probab.}, 30(3):1091--1129, 06 2020.

\bibitem[BFHM20]{benjamini2020}
I. Benjamini, L.~R. Fontes, J. Hermon, and F.~P. Machado.
\newblock On an epidemic model on finite graphs.
\newblock {\em Ann. Appl. Probab.}, 30(1):208--258, 02 2020.

\bibitem[{Big}76]{Big76}
J.~D. {Biggins}.
\newblock {The first- and last-birth problems for a multitype age-dependent
  branching process.}
\newblock {\em {Adv. Appl. Probab.}}, 8:446--459, 1976.

\bibitem[Big95]{Big95}
J.~D. Biggins.
\newblock The growth and spread of the general branching random walk.
\newblock {\em Ann. Appl. Probab.}, 5(4):1008--1024, 1995.

\bibitem[Big97]{Big97}
J.~D. Biggins.
\newblock How fast does a general branching random walk spread?
\newblock In {\em Classical and modern branching processes ({M}inneapolis,
  {MN}, 1994)}, volume~84 of {\em IMA Vol. Math. Appl.}, pages 19--39.
  Springer, New York, 1997.

\bibitem[Big10]{Big_branching_out}
J.~D. Biggins.
\newblock Branching out.
\newblock In {\em Probability and Mathematical Genetics: Papers in Honour of
  Sir John Kingman}, pages 112--133. Cambridge University Press, 2010.

\bibitem[BK20]{frogL}
V. Bezborodov and T. Krueger.
\newblock Linear and superlinear spread for continuous-time frog model.
\newblock arXiv:2008.10585, 2020.

\bibitem[BM14]{BM14}
J. B\'{e}rard and P. Maillard.
\newblock The limiting process of {$N$}-particle branching random walk with
  polynomial tails.
\newblock {\em Electron. J. Probab.}, 19:no. 22, 17, 2014.

\bibitem[CD16]{CD16}
S. {Chatterjee} and P.~S. {Dey}.
\newblock {Multiple phase transitions in long-range first-passage percolation
  on square lattices}.
\newblock {\em {Commun. Pure Appl. Math.}}, 69(2):203--256, 2016.

\bibitem[Dei03]{Dei03}
M. Deijfen.
\newblock Asymptotic shape in a continuum growth model.
\newblock {\em Adv. in Appl. Probab.}, 35(2):303--318, 2003.

\bibitem[DGH{\etalchar{+}}18]{RecFrog18}
C. {D\"obler}, N. {Gantert}, T. {H\"ofelsauer}, S. {Popov}, and F. {Weidner}.
\newblock {Recurrence and transience of frogs with drift on \(\mathbb{Z}^d\).}
\newblock {\em {Electron. J. Probab.}}, 23:23, 2018.
\newblock Id/No 88.

\bibitem[DHL19]{DHL19}
M. {Deijfen}, T. {Hirscher}, and F. {Lopes}.
\newblock {Competing frogs on \({\mathbb Z}^d\).}
\newblock {\em {Electron. J. Probab.}}, 24:17, 2019.
\newblock Id/No 146.

\bibitem[Dur83]{Dur83}
R. Durrett.
\newblock Maxima of branching random walks.
\newblock {\em Z. Wahrsch. Verw. Gebiete}, 62(2):165--170, 1983.

\bibitem[FZ12]{FZ12}
M. Fang and O. Zeitouni.
\newblock Branching random walks in time inhomogeneous environments.
\newblock {\em Electron. J. Probab.}, 17:no. 67, 18, 2012.

\bibitem[Gan00]{Gan00}
N. Gantert.
\newblock The maximum of a branching random walk with semiexponential
  increments.
\newblock {\em Ann. Probab.}, 28(3):1219--1229, 2000.

\bibitem[GM08]{GM08}
J.-B. Gou{\'e}r{\'e} and R. Marchand.
\newblock Continuous first-passage percolation and continuous greedy paths
  model: linear growth.
\newblock {\em Ann. Appl. Probab.}, 18(6):2300--2319, 2008.

\bibitem[GNR17]{GNR17}
A. {Ghosh}, S. {Noren}, and A. {Roitershtein}.
\newblock {On the range of the transient frog model on \(\mathbb{Z}\).}
\newblock {\em {Adv. Appl. Probab.}}, 49(2):327--343, 2017.

\bibitem[HJJ16]{HJJ16}
C. {Hoffman}, T. {Johnson}, and M. {Junge}.
\newblock {From transience to recurrence with Poisson tree frogs.}
\newblock {\em {Ann. Appl. Probab.}}, 26(3):1620--1635, 2016.

\bibitem[HJJ17]{RecFrog17}
C. {Hoffman}, T. {Johnson}, and M. {Junge}.
\newblock {Recurrence and transience for the frog model on trees.}
\newblock {\em {Ann. Probab.}}, 45(5):2826--2854, 2017.

\bibitem[{Jun}20]{Junge20}
M. {Junge}.
\newblock {Critical percolation and \(\mathrm{A + B \rightarrow 2A}\)
  dynamics}.
\newblock {\em {J. Stat. Phys.}}, 181(2):738--751, 2020.

\bibitem[KRS12]{KRS12}
H. {Kesten}, A.~F. {Ram\'{\i}rez}, and V. {Sidoravicius}.
\newblock {Asymptotic shape and propagation of fronts for growth models in
  dynamic random environment.}
\newblock In {\em {Probability in complex physical systems. In honour of Erwin
  Bolthausen and J\"urgen G\"artner. Selected papers based on the presentations
  at the two 2010 workshops}}, pages 195--223. Berlin: Springer, 2012.

\bibitem[KS05]{KesS05}
H. Kesten and V. Sidoravicius.
\newblock The spread of a rumor or infection in a moving population.
\newblock {\em Ann. Probab.}, 33(6):2402--2462, 2005.

\bibitem[KS08]{KesS08}
H. Kesten and V. Sidoravicius.
\newblock A shape theorem for the spread of an infection.
\newblock {\em Ann. of Math. (2)}, 167(3):701--766, 2008.

\bibitem[Mal15]{Mal15}
B. Mallein.
\newblock Maximal displacement in a branching random walk through interfaces.
\newblock {\em Electron. J. Probab.}, 20:no. 68, 40, 2015.

\bibitem[{Ros}17]{Brownian_Frog_17}
J. {Rosenberg}.
\newblock {The frog model with drift on \(\mathbb{R} \).}
\newblock {\em {Electron. Commun. Probab.}}, 22:14, 2017.
\newblock Id/No 30.

\bibitem[RS04]{stocCombust}
A.~F. {Ram\'{\i}rez} and V. {Sidoravicius}.
\newblock {Asymptotic behavior of a stochastic combustion growth process.}
\newblock {\em {J. Eur. Math. Soc. (JEMS)}}, 6(3):293--334, 2004.

\bibitem[vdHK17]{ExplFpp}
R. van~der Hofstad and J. Komjathy.
\newblock Explosion and distances in scale-free percolation.
\newblock 2017.
\newblock preprint; arXiv:1706.02597 [math.PR].

\end{thebibliography}

\end{document}